\documentclass[12pt]{article}
\usepackage{amssymb,latexsym,amsmath}
\usepackage{graphicx,times}
\usepackage{hyperref}
\newtheorem{The}{Theorem}[section]
\newtheorem{Def}{Definition}[section]
\newtheorem{lem}{Lemma}[section]

\numberwithin{equation}{section}
\begin{document}
\begin{center}
{\LARGE {\bf  Simple way to prove compactness of closed intervals in simply ordered set with order topology}}
\vskip 1cm
{\Large {\bf Sachin B. Bhalekar}}\\
{\large {\bf Department of Mathematics, \\
Shivaji University, \\
Vidyanagar, Kolhapur - 416004, India}}\\
\medskip
Email Address: sachin.math@yahoo.co.in, sbb\_maths@unishivaji.ac.in\\
\end{center}

\begin{abstract}
In this note, we present a simpler way to prove the compactness of the closed intervals in simply ordered set with order topology. 
\end{abstract}

\section{Introduction}
Compactness of a space is described using open sets only and it is preserved under continuous functions. Hence it is a topological property. These spaces are very important in Mathematical analysis (\cite{Mun}, Page 147). The historical developments and original motivation of compactness is discussed in a nice paper \cite{sun} by Sundstr\''om. This article deals with a simple proof of compactness of closed intervals in a simply ordered set with order topology. 
\section{Preliminaries}
In this section, we discuss some basic definitions and results from the literature \cite{Mun, Sim}.
\begin{Def}
Let $X$ be a set. A topology on $X$ is a collection $\tau$ of subsets of $X$ which is closed under arbitrary unions and finite intersections and contains an empty set $\phi$ and the set $X$. Elements of $\tau$ are called open sets in $X$.
\end{Def}

\begin{Def}
A set $X$ is called simply ordered if there is a relation $``<"$ on $X$ satisfying the following properties:\\
1. If $x, y\in X$ with $x\ne y$ then either $x<y$ or $y<x$.\\
2. The relation $x<x$ does not hold for any $x\in X$.\\
3. If $x<y$ and $y<z$ the $x<z$.
\end{Def}

\begin{Def}
If $X$ is simply ordered set with the order relation $<$ and if $a$ and $b$ are elements in $X$ with $a<b$ then the open interval $(a, b)$ in $X$ is a set $\left\{x\in X | a<x<b\right\}$ and closed interval $[a, b]$ in $X$ is a set $\left\{x\in X | a\leq x \leq b\right\}$. 
\end{Def}

\begin{Def}
An ordered set $X$ is said to have the least upper bound (lub) property if every nonempty subset of $X$ that is bounded above has lub. 
\end{Def}

\begin{Def}
Let $X$ be a simply ordered set with more than one elements. If $\mathcal{B}$ is collection of \\
1. all open intervals $(a, b)$ in $X$,\\
2. all intervals of the form $[a_0, b)$ (respectively, $(a, b_0]$), where $a_0$ (respectively, $b_0$) is the smallest (respectively, largest) element, if any, of $X$\\
then it is a basis for a topology on $X$, called the order topology.
\end{Def}

\begin{Def}
A sequence $\left\{x_n\right\}_{n=1}^{\infty}$ of points of the topological space $X$ is said to converge to the point $x$ of $X$ if for each open set $U$ containing $x$, $\exists$ $n_0\in \mathbb{Z}_+$ such that $x_n\in U,$ $\forall$ $n>n_0$.
\end{Def}

\begin{lem}
Let $Y$ be a subspace of topological space $X$. Then $Y$ is compact if and only if every covering of $Y$ by sets open in $X$ contains a finite subcollection covering $Y$.
\end{lem}

\begin{The}
A topological space is compact if every open cover by basis elements has a finite subcover.
\end{The}


\section{Main Result}
Heine-Borel theorem shows that every closed and bounded interval of real line is compact in standard topology. There is another way discussed in \cite{stack} which divides the interval in two parts and uses the nested intervals theorem. However, if we have arbitrary topological space with order topology then we don't have a metric on it. These approaches will not help us to prove the compactness of closed intervals in such spaces. In \cite{Mun}, the proof is given in four steps which is not easy. 
In this section, we present a new and simple way to prove this Theorem 27.1 in \cite{Mun}.

\begin{The} \cite{Mun}
Let $X$ be a simply ordered set having least upper bound property. In the order topology, each closed interval in X is compact. 
\end{The}
\textbf{Proof:} 
Suppose that the interval $[a, b]\in X$ is covered by the collection of basis elements $\left\{I_\alpha\right\}_{\alpha\in J}$ in $X$, where $I_\alpha=\left(a_\alpha, b_\alpha\right)$ and $J$ is an index set.
\par Since $a, b\in [a, b]\subset \cup_{\alpha\in J}I_\alpha$, there exist indices $\alpha_1$ and $\alpha_2$ in $J$ such that $a\in I_{\alpha_1}$ and $b\in I_{\alpha_2}$. 
If $I_{\alpha_1}\cap I_{\alpha_2}$ is not empty then $[a, b]\subseteq I_{\alpha_1} \cup I_{\alpha_2}$ and we are done. If not, then as shown in Figure 1, $[a, b]\subseteq I_{\alpha_1} \cup I_{\alpha_2}\cup \left[b_{\alpha_1}, a_{\alpha_2}\right]$.

\begin{figure}[h]
\centering
\includegraphics[width=0.7\linewidth]{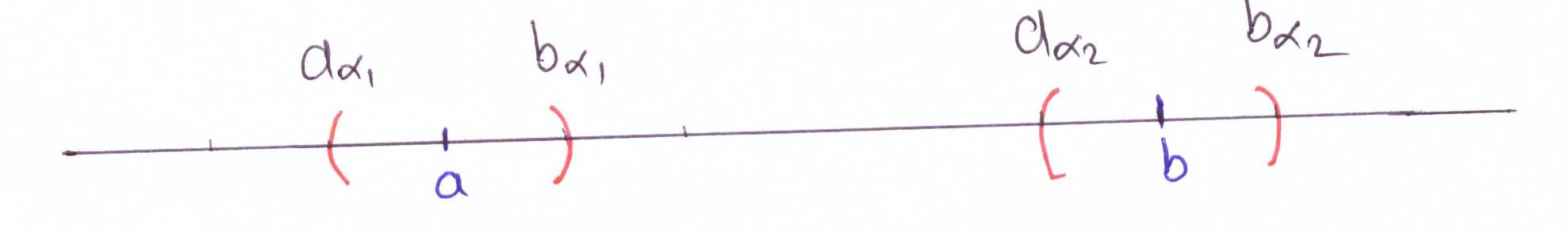}
\caption{$[a, b]\subseteq I_{\alpha_1} \cup I_{\alpha_2}\cup \left[b_{\alpha_1}, a_{\alpha_2}\right]$}
\label{fig:001}
\end{figure}

\par Now, consider the interval $\left[b_{\alpha_1}, a_{\alpha_2}\right]$ which is proper subset of $[a, b]$. By similar arguments, there exist indices $\alpha_3$ and $\alpha_4$ in $J$ such that $b_{\alpha_1}\in I_{\alpha_3}$ and $a_{\alpha_2}\in I_{\alpha_4}$.
If $I_{\alpha_3}\cap I_{\alpha_4}$ is not empty then $[a, b]\subseteq  \cup_{j=1}^{4} I_{\alpha_j}$ and we are done. If not, then $[a, b]\subseteq \cup_{j=1}^{4} I_{\alpha_j}\cup \left[b_{\alpha_3}, a_{\alpha_4}\right]$.
\par Continuing this procedure, we construct a bounded sequence $b_{\alpha_1}< b_{\alpha_3} < b_{\alpha_5}< \cdots$ in $[a, b]$. Since $X$ has least upper bound (lub) property, there exists lub $b_*$ of the set $\left\{b_{\alpha_1}, b_{\alpha_3}, b_{\alpha_5}, \cdots \right\}$. Note that 
\begin{equation*}
\lim_{j\longrightarrow \infty}b_{\alpha_{2j-1}}=b_* < b.
\end{equation*}
Since $b_*\in [a, b]$, there exists $\alpha_* \in J$ such that $b_*\in I_{\alpha_*}$. By definition of convergence of sequence, $\exists\, n_0$ such that 
\begin{equation*}
b_{\alpha_{2j-1}} \in I_{\alpha_*}, \quad  \forall j>n_0.
\end{equation*}
Therefore,
\begin{equation*}
[a, b]\subseteq \cup_{j=1}^{2n_0+1} I_{\alpha_j}\cup I_{\alpha_*}.
\end{equation*}
Thus, $[a, b]$ is covered by finitely many basis elements. Hence, $[a, b]$ is compact set in $X$.

\section{Conclusion}
In this article, we provided a simple proof of compactness of closed intervals in simply order set with order topology. We hope that the students will find this alternate way simpler than the one available in the literature.

\end{document}